\documentclass{ifacconf}

\usepackage{graphicx}      
\usepackage{natbib}        
\usepackage{amsmath}
\usepackage{amsfonts}

\usepackage{mathrsfs}  
\usepackage{physics}
\usepackage{mathtools}
\mathtoolsset{showonlyrefs}

\newtheorem{ass}{Assumption}[section]

\usepackage{xcolor}

\begin{document}
\begin{frontmatter}
\vspace{-1em}
\title{PINNs and GaLS: A Priori Error Estimates for Shallow Physics Informed Neural Networks Applied to Elliptic Problems} 
\vspace{-0.5em}
\author[First]{U. Zerbinati}
\vspace{-0.5em}
\address[First]{King Abdullah University Of Science and Technology, 
   Thuwal, SA (e-mail: umberto.zerbinati@kaust.edu.sa).}

\begin{abstract}                
Physics Informed Neural Networks (PINNs) have recently gained popularity for solving partial differential equations, given the fact they escape the curse of dimensionality.
In this paper, we present Physics Informed Neural Networks as an underdetermined point matching collocation method then expose the connection between Galerkin Least Squares (GALS) and PINNs, to develop an a priori error estimate, in the context of elliptic problems. In particular, techniques that belong to the realm of least squares finite elements and Rademacher complexity analysis are used to obtain the error estimate.
\end{abstract}

\end{frontmatter}

\section{Introduction}
\vspace{-1em}
A number of algorithms have been recently developed to solve partial differential equations by using deep and shallow neural networks, in order to escape the curse of dimensionality, e.g. \cite{MolinaroRadiative, Nature, WeinanCOD, Selectnet, NatureCom, WeinanNonLinear, WeinanBarronMeanField}. 
In particular, three methods are dominating the landscape of numerical solutions of partial differential equations (PDEs) by means of deep learning: the deep Ritz method more information in \cite{WeinanDeepRitz} and \cite{WeinanNonLinear}; the finite neuron method in \cite{FiniteNeuron} and \cite{JinchaoRelu}; and PINNs in \cite{DeepXDE} and \cite{Nature}.
Key differences among the three methods are that the first two are based on the variational formulation of the PDE, while the third one is based on the strong formulation of the problem.  Physics Informed Neural Networks (PINNs) can be be seen as an underdetermined point matching collocation method, as described in \cite{LSQBochev153} and \cite{Bochev}, Section 12.4. In fact the idea underlying PINNs is to minimize the least square residual evaluated at a fixed number of collocation points.
Obtaining a discrete least square principle by first discretizing the integral in the energy functional using a quadrature rule is not a new concept. It can be found for example in \cite{Bochev230}. The key difference between PINNs and a standard least square collocation method is that with the latter, one tends to choose a greater number of collocation points compared to the number of basis function in the discrete space, which is not the case for PINNs. 
Usually a priori error estimates are developed using Rademacher complexity analysis, as shown in \cite{APrioriJinchao} and \cite{KarniadakisPINNs}, however this technique results in underestimating the convergence rate of the method.
 In order to deal with this problem another approach has been developed, which is more similar to an a posteriori error estimate. This idea is presented in \cite{MolinaroEstimates} and \cite{MultiLevel}.
The aim of this article is to expose the connection between Galerkin Least Squares (GaLS) and PINNs, and to develop an a priori error estimate using techniques that belong to the realm of the least square finite elements for PINNs. To achieve this, we give a brief overview on GaLS in section 2, essentially following \cite{Bochev}.
In section 3, we present the notion of PINNs and present some results that will be used in section 4. These concepts can be found in \cite{DeepXDE} and \cite{APrioriJinchao}.
In section 4, we develop an a priori error estimate for PINNs when dealing with an elliptic operator, building on the idea that were introduced by \cite{Cai}. The novelty of this paper is that it applies ideas used to analyse the convergence of GaLS least squares methods to PINNs.
\vspace{-0.4 em}
\section{Galerkin Least Squares}
\vspace{-0.8em}
Let $X$ and $Y$ be two Hilbert spaces and consider a Fredholm operator $Q\in\mathcal{L}(X,Y)$, we will now focus our attention on the following problem,
\begin{equation}
find \; u\in X \; such \; that \; Qu=F. \label{eq:P1}
\end{equation}
\noindent
We will work under the following assumption to simplify the discussion, $\dim \Big(N(Q)\Big)=0$, where $N(Q)$ is the null space of $Q$.
We will focus our interest on the \textbf{residual energy functional} and we consider the following minimization principle,
\begin{equation}
	J(u;F) = \norm{Qu-F}^2_Y,
	\qquad
	u = \underset{v\in X}{\arg\min} \; J(v;F) \label{eq:M1}.
\end{equation}
The above minimization principle admits a unique minimizer. To show this we use the following well known result.
\begin{thm}\label{thm:ConvexMin}
	Given a reflexive Banach space $X$ and a continuous and strongly convex functional, $J: X \to \mathbb{R},$
	if the following conditions are satisfied,
	\begin{enumerate}
		\item $\underset{\norm{x}\to \infty}{lim}J(x)=\infty$,
		\item $K$ is a closed convex subset of $X$,
	\end{enumerate}
	than it exists a unique element $x^*\in K$ such that,
	\begin{equation}
	J(x^*)=\underset{x\in K}{\inf}J(x).
	\end{equation}
\end{thm}
\begin{thm}
	The minimization principle \eqref{eq:M1} has a unique minimizer $u^*\in X$ for any $F\in Y$.
\end{thm}
\begin{pf}
	First we want to prove the ellipticity of the operator $Q$. To do this we notice that by assumption the range of $Q$, denoted $R(Q)$, is closed in $Y$ and therefore $\Big(R(Q),(\cdot,\cdot)_Y\Big)$ is a Hilbert space. Now we consider the restricted mapping,
	\begin{equation}
	Q:X\to R(Q),
	\end{equation}
	such mapping is a bijection because of the assumption $\dim \Big(N(Q)\Big)=0$ and the fact we have taken as codomain of the operator $R(Q)$. Since the mapping is a bijection one can use the bounded inverse theorem to state that $Q^{-1}:R(Q)\to X$ is a linear bounded operator and therefore,
	\begin{equation}
	\norm{Qu}_Y\geq C_1 \norm{u}_X.\label{eq:Ellipiticity}
	\end{equation}
	Now we proceed to prove the coercivity for $J(u;F)$, which is property 1 in \textit{Theorem \ref{thm:ConvexMin}} statement.
	\begin{gather}
	J(u;F) = \Big(Qu-F,Qu-F\Big)_Y=\\(Qu,Qu)_Y-2(Qu,F)_Y+(F,F)_Y=\\\norm{Qu}^2_Y+\norm{F}^2_Y-2(Qu,F)_Y
	\end{gather}
	using Young inequality with $\varepsilon=2$ and \eqref{eq:Ellipiticity} we get the following inequality,
	\begin{gather}
	J(u;F)\geq \frac{1}{2}\norm{Qu}_Y^2 - \norm{F}^2_Y\geq \frac{C_1^2}{2} \norm{u}^2_X- \norm{F}^2_Y,
	\end{gather}
	which implies $\underset{\norm{u}\to \infty}{\lim}J(x)=\infty$. The only thing left  to prove is that the functional $J(u;F)$ is strictly convex, to do this we consider $u,v\in X$ and $t\in [0,1]$ and we evaluate,
	\begin{gather}
	j(t) = J(tu+(1-t)v;F) =\\ t^2 (Qu-F,Qu-F)_Y+(1-t)^2 (Qv-F,Qv-F)_Y+\\t(1-t)(Qu-F,Qv-F)_Y.
	\end{gather}
	If the coefficient of the quadratic term is greater than zero then the functional $J(u;F)$ is convex and this is the case because,
	\begin{gather}
	(Qu-F,Qu-F)_Y-2(Qu-F,Qv-F)_Y\\+(Qv-F,Qv-F)_Y\geq 0
	\end{gather}
	where we have used Young's inequality to obtain the last bound. In particular we notice that if 
	$$(Qu-F,Qu-F)_Y-2(Qu-F,Qv-F)_Y+(Qv-F,Qv-F)_Y$$
	is null then,
	\begin{gather}
	(Qu-Qv,Qu-Qv)=0 \Leftrightarrow Q(u-v)=0 \Leftrightarrow u=v,
	\end{gather}
	where the last implication comes from the fact that we assume $\dim \Big(N(Q)\Big)=0$. To conclude one needs to apply the \textit{Theorem \ref{thm:ConvexMin}}.
\end{pf}
\noindent
In particular one can characterize the minimizer of \eqref{eq:M1} using the following proposition.
\begin{prop}
	The minimizer $u^*$ of \eqref{eq:M1} solves the following variational equation,
	\begin{equation}
	find\; u \in X \; such \; that \; a(u,v)=G(v)\; \forall v \in X,\label{eq:V1}
	\end{equation}
	where $a(u,v)$ is the scalar product $\Big(Qu,Qv\Big)_Y$ and $G(v)$ is the following linear functional
	\begin{gather}
	G:X\to Y\\
	v\mapsto (F,Qv)_Y.
	\end{gather}
\end{prop}
\begin{pf}
	We define the following function of a real variable $j(t) = J(u^*+tv;F)$ and we notice that since $u^*$ is a minimum for $J(u;F)$ then $j'(0)=0$. Expanding $j'(t)$ we get,
	\begin{gather}
	j'(t)=2t(Qv,Qv)_Y+2(Qv,Qu-F)_Y.
	\end{gather}
	Imposing $j'(0)=0$  for all $v \in X$ one finds \eqref{eq:V1}.
\end{pf}
\noindent
Now that we have shown that the continuous minimization principle \eqref{eq:M1} has a unique minimizer and such a minimizer is characterized by solving the variational equation \eqref{eq:V1} we focus our attention on the discrete problem associated with \eqref{eq:M1}. Usually in the context of GaLS one chooses a conforming finite dimensional subspace $X_h\subset X$ and consider as discrete least square principle \eqref{eq:M1} the energy functional evaluated on $X_h$, i.e.
\begin{equation}
J(u;F) = \norm{Qu-F}^2_Y,
\qquad
u_h = \underset{v\in X_h}{\arg\min} \; J(v;F) \label{eq:M1h}.
\end{equation}
The advantage of the above mentioned approach is that to produce an a priori error estimate one can follow the standard practice of using Cea's lemma and interpolation error estimate in $X_h$.
Another approach would be to use directly a discrete minimization principle, i.e.
\begin{equation}
J^h(u;F) = \norm{Qu-F}^2_h,
\qquad
u_h = \underset{v\in X_h}{\arg\min} \; J^h(v;F) \label{eq:M2}.
\end{equation}
where $\norm{\cdot}_h$ is the discrete energy norm associated with a discrete inner product,
$(\cdot,\cdot)_h:X_h\times X_h \to \mathbb{R}.$
A useful assumption to develop a priori error estimate for the above described discrete energy functional is that the discrete energy functional can be extended to $u\in X$. We will further assume that two positive semi-definite bilinear form $e(\cdot,\cdot)$ and $\epsilon(\cdot,\cdot)$ exist, such that,
\begin{gather}
J^h(u^h,F)=\frac{1}{2} \Bigg((u_h,u_h)_h+(u,u)_h+\epsilon(u,u)\Bigg)-(u,u_h)_h\\-e(u,u_h)\qquad \forall u\in X, \; \forall u_h \in X_h,\label{eq:DiscreteEnergyFunctional}
\end{gather}
under this assumption and proceed to prove an a priori error estimate for our discrete least square principle \eqref{eq:M2}. It is important to notice that we can proceed as we have done at the beginning of this section in order to prove the existence of a minimizer for \eqref{eq:M2}.
\begin{thm}\label{lem:Bochev}
	Under the above assumption and letting $u\in X$ be the solution of \eqref{eq:M1}, $u_h\in X_h$ be the unique solution of \eqref{eq:M2}, then the following error estimate holds,
	\begin{equation}
	\norm{u-u_h}_h \leq \underset{v \in X_h}{\inf} \norm{u-v}_h+\underset{v\in X_h}{\sup}\frac{e(u,v)}{\norm{v}_h}.\label{eq:Bochev}
	\end{equation}
	In particular if the least square functional $J^h$ is r-consistent, i.e. it exist $r>0$ such that,
	\begin{equation}
	\underset{v\in X_h}{\sup} \frac{e(u,v)}{\norm{v}_h}\leq C(u)h^r,
	\end{equation}
	then the above error estimate becomes,
	\begin{equation}
	\norm{u-u_h}_h \leq \underset{v \in X_h}{\inf} \norm{u-v}_h + C(u)h^r.
	\end{equation}
\end{thm}
\begin{pf}
	Let $u_h^{\bot}$ be the orthogonal projection of $u$ onto $X^h$ with respect to the discrete scalar product $(\cdot,\cdot)_h$.
	We know that the following chain of inequalities holds,
	\begin{gather}
	\norm{u_h-u}_h\leq\norm{u_h^{\bot}-u}_h+\norm{u_h-u_h^{\bot}}_h\\\leq \underset{v \in X_h}{\inf} \norm{u-v}_h+\norm{u_h-u_h^{\bot}}_h.
	\end{gather} 
	Therefore we are left to estimate the term $\norm{u_h-u_h^{\bot}}_h$. In order to achieve such estimate we observe that,
	\begin{equation}
	(u_h-u,v)_h = e(u,v) \qquad \forall v \in X_h
	\end{equation}
	is the variational equation associated with minimizing \eqref{eq:DiscreteEnergyFunctional} and given the fact that $u_h^{\bot}$ is the orthogonal projection of $u$ one has,
	\begin{equation}
	(u_h^{\bot}-u,v)_h=0 \qquad \forall v \in X_h.
	\end{equation}
	Subtracting the last two equations and using the definition of norm induced by the linear operator we get,
	\begin{gather}
	(u_h^{\bot}-u_h,v)_h = e(u,v)\qquad v \in X_h\\
	\norm{u_h^{\bot}-u_h}_h=\underset{v\in X_h}{\sup} \frac{(u_h^{\bot}-u_h,v)_h}{\norm{v}_h}=\underset{v\in X_h}{\sup} \frac{e(u,v)_h}{\norm{v}_h}.
	\end{gather}
\end{pf}
\noindent
We will split the contribution of the differential operator on the inside of the domain $\Omega$ and on its boundary $\partial\Omega$ obtaining the following version of \eqref{eq:P1},
\begin{gather}
\label{eq:SplitedOperator}
Q:X\to Y := A(\Omega)\times B(\partial\Omega)\\
u\mapsto \Big(\mathcal{L}u,\mathcal{B}u\Big),
\end{gather}
where $A(\Omega)$ and $B(\partial\Omega)$ are Hilbert spaces representing respectively component of $Y$ on $\Omega$ and $\partial\Omega$.

\section{Physics Informed Neural Network}
\vspace{-0.5em}
In this section we will introduce the notion of Physics Informed Neural Networks (PINNs) and the correct functional setting where to study PINNs.
\begin{defn}[Forward Neural Network]
	We say $u_h^{L}(\vec{x}):\mathbb{R}^{d_{in}}\to \mathbb{R}^{d_{out}}$ is a $(L-1)$ hidden layer forward neural network (FNN) with $N_\ell$ neurons in the $\ell-th$ layer and activation function $\sigma:\mathbb{R}\to\mathbb{R}$, if the action $u_h^{L}(\vec{x})$ is recursively defined as,
	\begin{enumerate}
		\item $u_h^0(\vec{x})= \vec{x}\in \mathbb{R}^{d_{in}}$,
		\item $\forall 1 \leq \ell \leq L-1\qquad u_h^\ell(\vec{x}) = \sigma(W^\ell N^{\ell-1}(\vec{x})+\vec{b_\ell})\in \mathbb{R}^{N_\ell}$,
		\item $u_h^{L}(\vec{x}) = W^Lu_h^{L-1}(\vec{x})+\vec{b}^{L}\in \mathbb{R}^{d_{out}}$,
	\end{enumerate}
	where $W^\ell \in \mathbb{R}^{N_\ell\times N_{\ell-1}}$ and $\vec{b_\ell}\in\mathbb{R}^{N_\ell}$. In particular the value of the matrix $W^\ell$ will be called kernel parameters and the value of the vector $\vec{b}^\ell$ will be called bias parameters for the layer $\ell$. The parameters of the network will be denoted as $\vec{\theta}=\Big(\{W_\ell\}_{\ell=1}^L,\{\vec{b}_\ell\}_{\ell=1}^L\Big)$ and we indicate the dependence of the FNN on a particular choice of parameters $\hat{\vec{\theta}}$ writing $u^L_{h,\hat{\theta}}$. We will often omit both the dimension of the network and the dependence on the parameters writing simply $u_h$.
\end{defn}
\noindent
We will focus our attention on the logistic sigmoid, i.e. $\sigma(x) = \frac{1}{1+e^{-x}}$, and the hyperbolic tangent as activation functions. When we speak about {Physics Informed Neural Networks (PINNs), we mean a FNN trained to solve \eqref{eq:P1} using the following algorithm:
\begin{enumerate}
	\item Given a probability distribution of points on a set $A$, called $\mathcal{D}(A)$, we consider two vectors of realisations, i.e. 
	\begin{gather}\omega_i \sim \mathcal{D}(\Omega)\qquad \forall 1\leq i \leq N_\Omega \\
	  \beta_i \sim \mathcal{D}(\partial\Omega)\qquad \forall 1\leq i \leq N_{\partial\Omega}.
	\end{gather}
	\item We consider the following loss function, i.e.
	\begin{gather}
	\mathscr{L}(\theta;\omega_i,\beta_i)=\frac{1}{N_{\Omega}^{\gamma_1}}\sum_{i=1}^{N_\Omega} \norm{\mathcal{L}u^L_h({\omega_i})-f(\omega_i)}_a^2+\\\frac{1}{N_{\partial\Omega}^{\gamma_2}}\sum_{i=1}^{N_{\partial\Omega}} \norm{\mathcal{B}u^L_h({\beta_i})-g(\beta_i)}_b^2,\label{eq:LossFunction}
	\end{gather}
	where $\gamma_i$ are parameters depending on the dimension of $\Omega$, and $F=(f,\,g)\in A(\Omega)\times B(\partial\Omega)$.
	\item We train the FNN in order to find $\theta^*$ that minimize the above loss function.
\end{enumerate}
It is worth noticing that the optimization problem in step 3 of the PINNs algorithm is highly non convex with respect to $\vec{\theta}$, \cite{NP}, and the state of the art for training PINNs is to use a first-order descent method such as Adam or L-BFGS. 

In order to develop an a priori error estimate for PINNs we need a class of sufficient regularity for the solution of our PDE, thus we will here give a very brief brush on concepts that are treated in more detail in \cite{JinchaoFunctional},\cite{WeinanBarron},\cite{WeinanBarronMeanField} and \cite{SiegalVariation}.
Let $X$ be the Hilbert space that we have introduced in the previous section and $\mathbb{D}\subset X$ a dictionary of elements in $X$ such that,
$K_{\mathbb{D}}:= \underset{d\in \mathbb{D}}{\sup}\;\norm{d}_X<\infty.$\\
In the above setting one can define the functional space of interest for our work, as:
\begin{defn}[Barron Space]
	Let us consider the closure of the convex symmetric hull of $\mathbb{D}$,
	\begin{equation}
	B_1(\mathbb{D}) = \overline{\Big\{\sum_{j=1}^{n} a_jd_j \; : \; n\in\mathbb{N},d_j\in \mathbb{D}, \; \norm{\{a_j\}_{j=1}^n}_{\ell_1}\leq 1\Big\}}.
	\end{equation} 
	we define the Barron space and the associated norm as follow:
	\begin{gather}
	\norm{\cdot}_{\mathcal{K}(\mathbb{D})} = \inf\Big\{c>0 \; : \; cf\in B_1(\mathbb{D})\Big\},\\ \mathcal{K}(\mathbb{D})=\Big\{f \in X \; : \; \norm{f}_{\mathcal{K}(\mathbb{D})}<\infty\Big\}.
	\end{gather}
	In particular given the fact that $X$ is a Hilbert space it is obvious that $\mathcal{K}_{\mathbb{D}}$ is a subset of $X$.
\end{defn}
\noindent
The importance of the Barron space is that given a function in $\mathcal{K}(\mathbb{D})$ we have an a priori approximation estimates for shallow neural network. To formalize such a concept we notice that the set of shallow FNN, i.e. FNN with one hidden layer and null bias on the output layer, can be rewritten as
\begin{equation}
\Sigma_{N,M}(\mathbb{D})=\Big\{\sum_{j=1}^{N} a_jd_j \; : d_j\in \mathbb{D}\; and \; \norm{\{a_j\}_{j=1}^N}_{\ell_1}\leq M\Big\},
\end{equation}
where we imposed the additional constraint on the $\ell_1$ norm of the coefficients $a_j$ and chosen a particular dictionary,
$$\mathbb{D}_\sigma=\Big\{\sigma(\vec{w}_i\cdot x+\vec{b}_i)\; : \; \vec{w}_i\in \mathbb{R}^{d}\; and \; \vec{b}_i\in \mathbb{R}\Big\}.$$
Now since $X$ is a Hilbert space it is also a type 2 Banach space, we have the following approximation estimate.
\begin{thm}[\cite{Mauery}, \cite{DeVoreActaNumerica}]
	\label{thm:Maurey}
	$\;\;\;$
	Let $X$ be a Hilbert space and given $f\in\mathcal{K}(\mathbb{D_\sigma})$ we have the following approximation estimate:
	\begin{equation}
	\underset{f_N \in \Sigma_{N,M}(\mathbb{D}_\sigma)}{\inf}\norm{f-f_N}_X \leq C_{X}K_{\mathbb{D}_\sigma}\norm{f}_{\mathcal{K}_1(\mathbb{D}_\sigma)}N^{-\frac{1}{2}},
	\end{equation}
	where $C_X$ is the type-2 constant for the space $X$.
	Furthermore when $\sigma$ is a bounded activation function then $\mathbb{D_\sigma}$ is uniformly bounded in $\mathcal{K}(\mathbb{D_\sigma})$.
\end{thm}
\noindent
The above error estimate can be improved for certain types of activation functions. More details can be found in \cite{JinchaoFunctional}, \cite{Makovoz} and \cite{PinkusActa}.\\
In the next section we will use the concept developed in the previous section together with \textit{Theorem  \ref{thm:Maurey}} to prove an a priori error estimate, similarly to what is done using Céa's lemma in the context of classical conforming finite element methods.
To conclude this section we notice that a regularity theory for elliptic problems in Barron space has yet to be fully developed, in particular the sufficient condition to impose on the source function in order for the solution to live in $\mathcal{K}(\mathbb{D_\sigma})$ are not known, more information on this problem can be found in \cite{Representation} and \cite{WeinanBarronBruna}.

\vspace{-0.5em}
\section{A Priori Error Estimate}
\vspace{-0.85em}

In this section the a priori error estimate that is the aim of this paper will be presented.
In particular we will consider the following energy functional associated with the strong formulation \eqref{eq:SplitedOperator},
\begin{gather}
J(u;F) = \norm{Qu-F}_Y^2 = \int_Y \norm{Qu-F}_y \;d\vec{x}= \\\int_\Omega \norm{\mathcal{L}u-f}_a \;d\vec{x} + \int_{\partial\Omega} \norm{\mathcal{B}u-g}_b\;d\vec{x}\label{eq:EnergyFunc}
\end{gather}
where $Q=(\mathcal{L},\mathcal{B})$, the space $Y$ is the product space $A(\Omega)\times B(\partial\Omega)$ and the space $A$ and $B$ are Hilbert spaces equipped with the following scalar product and induced norm,
\begin{gather}
(u,v)_{A(\Omega)} = \int_{\Omega} \Big(u(\vec{x}),v(\vec{x})\Big)_a \; d\vec{x} \\ \norm{u}_{A(\Omega)}^2 = \int_{\Omega} \Big(u(\vec{x}),u(\vec{x})\Big)_a \; d\vec{x}=\int_{\Omega} \abs{u(\vec{x})}_a^2 \; d\vec{x},\\
(u,v)_{B(\partial\Omega)} = \int_{\partial\Omega} \Big(u(\vec{x}),v(\vec{x})\Big)_b \; d\vec{x} \\ \norm{u}_{B(\partial\Omega)}^2 = \int_{\partial\Omega} \Big(u(\vec{x}),u(\vec{x})\Big)_b \; d\vec{x}=\int_{\partial\Omega} \abs{u(\vec{x})}_b^2 \; d\vec{x}
\end{gather}
The discrete scalar product that gives us the discrete energy functional \eqref{eq:LossFunction} is:
\begin{gather}
(\cdot,\cdot)_h: X_h\times X_h\to \mathbb{R}\\
(u_h,v_h)_h\mapsto \frac{1}{N_{\Omega}^{\gamma_1}}\sum_{i=1}^{N_\Omega} \Big(u_h({\vec{\omega}_i}),v_h({\vec{\omega}_i})\Big)_a+\\ \frac{1}{N_{\partial\Omega}^{\gamma_2}}\sum_{i=1}^{N_{\partial\Omega}} \Big(u_h({\vec{\beta}_i}),v_h({\vec{\beta}_i})\Big)_b.\label{eq:DiscreteScalarProduct}
\end{gather}
In particular we notice that the above mapping is bilinear thanks to the fact that $(\cdot,\cdot)_a$ and $(\cdot,\cdot)_b$ are scalar products. Furthermore since $(\cdot,\cdot)_a$ and $(\cdot,\cdot)_b$ are scalar products one has that $(u_h,u_h)_h \geq 0$ for all $u_h\in X_h$. Last we have to prove that $(u_h,u_h)_h=0 \Leftrightarrow u_h\equiv 0$ in order to do this we will work under the following assumption.
\begin{ass}\label{ass:Assumption}
	Fixed $u_h\in X_h$ it exists a set of points $\Big\{\delta_1,\dots,\delta_{NDOF}\Big\}$ such that, $u_h(\delta_1)=\dots=u_h(\delta_{NDOF})=0 \Leftrightarrow u_h \equiv 0$.
	The above mentioned assumption can also be rephrased as $u_h(\delta_i)$ are unisolvent degrees of freedom for functions in $X_h$.
\end{ass}
\noindent
It's important to notice that in order for \textit{Assumption \ref{ass:Assumption}} to grant that $(u_h,u_h)_h=0 \Leftrightarrow u_h\equiv 0$ one needs to slightly modify the algorithm used to train PINNs, in particular at each training step we need to choose the points where we evaluate the loss function \eqref{eq:LossFunction} dynamically, in order for \textit{Assumption \ref{ass:Assumption}} to be verified. Furthermore \textit{Assumption \ref{ass:Assumption}} can't be verified for all activation functions. For example, when dealing with second order PDE if we consider a sigmoid activation function it is sufficient to choose points close enough to the value of the bias layer for \textit{Assumption \ref{ass:Assumption}} to be verified. However \textit{Assumption \ref{ass:Assumption}} can't be verified if we consider a RELU activation function.

We now proceed using Rademacher complexity tools as shown in \cite{APrioriJinchao} in order to bound the right hand side of \eqref{eq:Bochev}.
\begin{thm}[\cite{APrioriJinchao},\cite{Rademacher}]\label{thm:Quadrature}$\;\;$
	Given a class of functions $\mathscr{F}:\Omega\to \mathbb{R}$ and a collection of sample points $\{\omega_i\}_{i=1}^{N_\Omega}$,
	\begin{gather}
	\mathbb{E}_{\omega_i\sim \mathcal{D}(\Omega)} \underset{h\in \mathscr{F}}{\sup} \abs{\sum^{N_\Omega}_{i=1}\frac{h(\omega_i)}{N_\Omega}-\int_\Omega h(x)\; d\mathcal{D}(\Omega)}\leq 2R_{N_\Omega}(\mathscr{F})\\:=\mathbb{E}_{\omega_i\sim \mathcal{D}(\Omega)}\mathbb{E}_{\xi_i}\Bigg[\underset{h\in \mathscr{F}}{\sup} \frac{1}{N_\Omega}\sum^{N_\Omega}_{i=1}\xi h(\omega_i)\Bigg],
	\end{gather}
	where $\xi_i$ are Rademacher random variables and $\omega_i$ are uniformly distributed points that includes the $\delta_i,\dots,\delta_{NDOF}$ from \textit{Assumption \ref{ass:Assumption}}.
\end{thm}
\begin{cor}
	Given a class of functions $\mathscr{F}:\partial\Omega\to \mathbb{R}$ and a collection of sample points $\{\beta_i\}_{i=1}^{N_{\partial\Omega}}$,
	\begin{gather}
	\mathbb{E}_{\beta_i\sim \mathcal{D}(\partial\Omega)} \underset{h\in \mathscr{F}}{\sup} \abs{\frac{1}{N_{\partial\Omega}}\sum^{N_{\partial\Omega}}_{i=1}h(\beta_i)-\int_{\partial\Omega} h(x)\; d\mathcal{D}(\partial\Omega)}\leq\\ 2R_{N_{\partial\Omega}}(\mathscr{F}):=\mathbb{E}_{\beta_i\sim \mathcal{D}(\partial\Omega)}\mathbb{E}_{\xi_i}\Bigg[\underset{h\in \mathscr{F}}{\sup} \frac{1}{N_{\partial\Omega}}\sum^{N_{\partial\Omega}}_{i=1}\xi h(\beta_i)\Bigg],
	\end{gather}
	where $\xi_i$ are Rademacher random variables and $\omega_i$ are uniformly distributed points that include the $\delta_i,\dots,\delta_{NDOF}$ from \textit{Assumption \ref{ass:Assumption}}.
\end{cor}
From now on when writing the expected value we will not specify the dependence on the distribution.
\begin{thm}[\cite{APrioriJinchao}]\label{thm:RademacherXu}
	Let $\mathscr{F},\mathscr{S}$ be classes of functions from $\Omega$ into $\mathbb{R}$, then the following statement holds,
	\begin{enumerate}
		\item $R_{N_p}(conv(\mathscr{F}))=R_{N_p}(\mathscr{F})$.
		\item $R_{N_p}(\mathscr{F}+\mathscr{S})=R_{N_p}(\mathscr{F})+R_{N_p}(\mathscr{S})$, where $$\mathscr{F}+\mathscr{S} = \Big\{h(x)+g(x)\;:\; h\in \mathscr{F},\; g \in \mathscr{S} \Big\}.$$
		\item $R_{N_p}(\Phi\circ\mathscr{F})\leq LR_{N_p}(\mathscr{F})$, where $$\Phi\circ \mathscr{F} = \Big\{\Phi\circ h(x)\;:\; h\in \mathscr{F}\Big \}$$ and $\Phi:\mathbb{R}\to\mathbb{R}$ is Lipschitz function with Lipschitz constant $L$.
		\item $R_{N_p}(f\cdot\mathscr{F})\leq \norm{f}_{\mathcal{L}^\infty(\Omega)}R_{N_p}(\mathscr{F})$, where $$f\cdot\mathscr{F} = \Big\{f(x)h(x)\;:\; h\in \mathscr{F}\Big \}$$ and $f:\Omega\to\mathbb{R}$ is a fixed function.
	\end{enumerate}
\end{thm}
\noindent
We will from now on assume that the differential operator $Q:X\to Y$ from \eqref{eq:P1} is an elliptic operator and we focus on the class of energy functional considered in \eqref{eq:EnergyFunc},
\begin{gather}
\mathcal{L}u := \sum_{\abs{\alpha}\leq m} a_\alpha(x)\partial^\alpha u, 
\qquad 
\mathcal{B}u := \sum_{\abs{\alpha}\leq m} b_\alpha(x)\partial^\alpha u,\label{eq:EllOp}
\\
\mathscr{F}_{N,M}^{\Omega}=\Bigg\{\Big(\mathcal{L}u-f,\mathcal{L}u-f\Big)_a \;:\; u \in \Sigma_{N,M}\Bigg\},
\\
\mathscr{F}_{N,M}^{\partial\Omega}=\Bigg\{\Big(\mathcal{B}u-g,\mathcal{B}u-g\Big)_b \;:\; u \in \Sigma_{N,M}\Bigg\}.
\end{gather}
\begin{thm}
	We will assume that $\norm{a_\alpha}_{\mathcal{L}^{\infty}(\Omega)}\leq K_a$, $\norm{b_\alpha}_{\mathcal{L}^{\infty}(\Omega)}\leq K_b$, $f\in\mathcal{L}^{\infty}(\Omega)$ and $g\in\mathcal{L}^{\infty}(\partial\Omega)$ then we have the following bound on the Rademacher complexity of $\mathscr{F}_{N,M}$,
	\begin{gather}
	R_{N_\Omega}(\mathscr{F}_{N,M}^{\Omega})\leq MK_aK_{\mathbb{D}}\sum_{\abs{\alpha}\leq m}R_{N_{\Omega}}(\partial^\alpha\mathbb{D})+\\2MK_a\norm{f}_{\mathcal{L}^\infty(\Omega)}\sum_{\abs{\alpha}\leq m}R_{N_{\Omega}}(\partial^\alpha\mathbb{D}),
	\\
	R_{N_{\partial\Omega}}(\mathscr{F}_{N,M}^{\partial\Omega})\leq MK_bK_{\mathbb{D}}\sum_{\abs{\alpha}\leq m}R_{N_{\partial\Omega}}(\partial^\alpha\mathbb{D})+\\2MK_b\norm{g}_{\mathcal{L}^\infty(\partial\Omega)}\sum_{\abs{\alpha}\leq m}R_{N_{\partial\Omega}}(\partial^\alpha\mathbb{D}),
	\end{gather}
	provided that $u\mapsto (u,u)_y$ is a locally Lipschitz function.
\end{thm}
\begin{pf}
	We begin considering $$\mathscr{F}_{N,M}=\Bigg\{\Big(Qu-F,Qu-F\Big)_y \;:\; u \in \Sigma_{N,M}\Bigg\}$$ and develop the scalar product in order to obtain,
	\begin{gather}
	\Big(\sum_{\abs{\alpha}\leq m} a_\alpha(x)\partial^\alpha u-f,\sum_{\abs{\alpha}\leq m} a_\alpha(x)\partial^\alpha u-f\Big)_a=\\ \Big(\sum_{\abs{\alpha}\leq m} a_\alpha(x)\partial^\alpha u,\sum_{\abs{\alpha}\leq m} a_\alpha(x)\partial^\alpha u\Big)_a
	\\-2\Big(\sum_{\abs{\alpha}\leq m}a_\alpha(x)\partial^\alpha u,\,f\Big)_a-(f,f)_a,\\		
	\Big(\sum_{\abs{\alpha}\leq m} b_\alpha(x)\partial^\alpha u-g,\sum_{\abs{\alpha}\leq m} b_\alpha(x)\partial^\alpha u-g\Big)_b=\\ \Big(\sum_{\abs{\alpha}\leq m} b_\alpha(x)\partial^\alpha u,\sum_{\abs{\alpha}\leq m} b_\alpha(x)\partial^\alpha u\Big)_b
	\\-2\Big(\sum_{\abs{\alpha}\leq m}b_\alpha(x)\partial^\alpha u,\,g\Big)_b-(g,g)_b.
	\end{gather}
	Now using \textit{Theorem \ref{thm:RademacherXu}} and the fact that $u\mapsto(u,u)_a,\;u\mapsto(u,u)_b$  are locally Lipschitz then one has the following bound on the Rademacher complexity of $\mathscr{F}_{N,M}$,
	\begin{gather}
	R_{N_\Omega}(\mathscr{F}_{N,M}^{\Omega})\leq MK_aK_{\mathbb{D}}\sum_{\abs{\alpha}\leq m}R_{N_{\Omega}}(\partial^\alpha\mathbb{D})+\\2MK_a\norm{f}_{\mathcal{L}^\infty(\Omega)}\sum_{\abs{\alpha}\leq m}R_{N_{\Omega}}(\partial^\alpha\mathbb{D}),
	\\
	R_{N_{\partial\Omega}}(\mathscr{F}_{N,M}^{\partial\Omega})\leq MK_bK_{\mathbb{D}}\sum_{\abs{\alpha}\leq m}R_{N_{\partial\Omega}}(\partial^\alpha\mathbb{D})+\\2MK_b\norm{g}_{\mathcal{L}^\infty(\partial\Omega)}\sum_{\abs{\alpha}\leq m}R_{N_{\partial\Omega}}(\partial^\alpha\mathbb{D}).
	\label{eq:RademacherBound}
	\end{gather}
	In order to obtain the above inequality we have also used the fact that the Rademacher complexity of a constant is null.
\end{pf}
\begin{thm}[\cite{APrioriJinchao}]\label{thm:RademacherBoundXu}
	Assuming the activation function $\sigma$ lives in $\mathcal{W}^{m+1,\infty}$, then for any $\alpha$ such that $\abs{\alpha}\leq m$ one has the following estimate for the Rademacher complexity,
	\begin{equation}
	R_{N_p}(\partial^\alpha \mathbb{D}_\sigma)\leq C{N_p}^{-\frac{1}{2}},
	\end{equation}
	where the constant $C$ does not depend on $N_p$.
\end{thm}
\begin{cor}\label{cor:Rademacher}
	In the hypothesis of \textit{Theorem \ref{thm:RademacherBoundXu}}, one has the following bound on the Rademacher complexity of $\mathscr{F}_{N,M}^{\Omega}$ and $\mathscr{F}_{N,M}^{\partial\Omega}$,
	\begin{gather}
	R_{N_\Omega}(\mathscr{F}_{N,M}^{\Omega})\leq MK_aK_{\mathbb{D}}N_{\Omega}^{-\frac{1}{2}}+2MK_a\norm{f}_{\mathcal{L}^\infty(\Omega)}N_{\Omega}^{-\frac{1}{2}},
	\\
	R_{N_{\partial\Omega}}(\mathscr{F}_{N,M}^{\partial\Omega})\leq MK_bK_{\mathbb{D}}N_{\partial\Omega}^{-\frac{1}{2}}+2MK_b\norm{g}_{\mathcal{L}^\infty(\partial\Omega)}N_{\partial\Omega}^{-\frac{1}{2}}.
	\end{gather}
\end{cor}
\noindent
Now we are ready to prove our a priori error estimate, because if we fix $\gamma_i=1$ for $i=1,2$ the above corollary together with \textit{Theorem \ref{thm:Quadrature}} tell us that our discrete energy functional is $\frac{1}{4}$-consistent and therefore one has from \textit{Theorem \ref{lem:Bochev}} that,
\begin{equation}
\mathbb{E}_{\omega,\beta}\norm{u-u_h}_h\leq\mathbb{E}_{\omega,\beta} \underset{v\in \Sigma_{N,M}}{\inf}\norm{u-v}_h+C_aN_\Omega^{-\frac{1}{4}}+C_bN_{\partial\Omega}^{-\frac{1}{4}},
\end{equation}
where $C_a$ is $\max\Big\{MK_aK_\mathbb{D_\sigma}, 2MK_a\norm{f}_{\mathcal{L}^\infty(\Omega)}\Big\}$ and $C_b$ is  $\max\Big\{MK_bK_\mathbb{D_\sigma}, 2MK_b\norm{g}_{\mathcal{L}^\infty(\partial\Omega)}\Big\}$.
Last we notice that once again using \textit{Theorem \ref{thm:Quadrature}} and \textit{Corollary \ref{cor:Rademacher}}, one has that $\mathbb{E}_{\omega,\beta}\norm{u-v}_h\leq \mathbb{E}_{\omega,\beta}\norm{u-v}_X+C_aN_\Omega^{-\frac{1}{4}}+C_bN_{\partial\Omega}^{-\frac{1}{4}}$ and therefore thanks to \textit{Theorem \ref{thm:Maurey}} the above equation becomes,
\begin{gather}
\mathbb{E}_{\omega,\beta}\norm{u-u_h}_h\leq C_{X}K_{\mathbb{D}_\sigma}\norm{u}_{\mathcal{K}(\mathbb{D}_\sigma)}N^{-\frac{1}{2}}\\+ 2C_aN_\Omega^{-\frac{1}{4}}+2C_bN_{\partial\Omega}^{-\frac{1}{4}}.
\end{gather}
where in this case $N$ is the number of neurons in the shallow layer.
\begin{thm}\label{thm:main}
	Considering \eqref{eq:P1} with $Q$ as in \eqref{eq:EllOp}, assuming $u\in \mathcal{K}(\mathbb{D_\sigma})$, $F\in\mathcal{L}^\infty(\Omega)$ and $\sigma\in\mathcal{W}^{m+1,\infty}$ then,
	\begin{gather}
	\mathbb{E}_{\omega,\beta}\norm{u-u_h}_X\leq 2C_{X}K_{\mathbb{D}_\sigma}\norm{u}_{\mathcal{K}_1(\mathbb{D}_\sigma)}N^{-\frac{1}{2}}\\+ 3C_aN_\Omega^{-\frac{1}{4}}+3C_bN_{\partial\Omega}^{-\frac{1}{4}},
	\end{gather}
	where in this case $N$ is the number of neurons in the shallow layer, $N_{\Omega}$ and $N_{\partial\Omega}$ are respectively the number of evaluation points on the boundary and inside of $\Omega$.
\end{thm}
\vspace{-1em}
\section{Conclusion}
\vspace{-0.8em}
A connection between Galerkin least square finite elements and physics informed neural network has been drawn and standard tools for least square finite elements methods have been used to develop an a priori error estimate for physics informed neural networks when solving \eqref{eq:P1}. A similar idea have been applied when considering the RELU activation function by \cite{Cai}. It is important to notice that the a priori error estimate presented here is suboptimal with respect to what is observed numerically when using a sigmoid activation function. This suboptimality can not be improved upon when using a sigmoid activation function, given the fact that in \cite{JinchaoFunctional} it has been proven that \textit{Theorem \ref{thm:Maurey}} is optimal for sigmoid activation functions. It is worth noting that suboptimality problems for a priori error estimates, seem to be common whenever Rademacher complexity is used and can be overcome using ideas presented in \cite{MolinaroEstimates}. It is important to notice that the a priori error estimate presented in \textit{Theorem \ref{thm:main}} only works if $u_h$ is the minimizer of the energy functional in $\mathcal{K}(\mathbb{D_\sigma})$. Given that PINNs training is a highly non convex optimization problem, there are no guarantees that one can find such a minimizer. On the other hand a priori error estimates presented in \cite{APrioriJinchao} depend on the use of a particular greedy algorithm but hold at each step of the minimization process. Last but not least numerical simulation showing the application of PINNs to different elliptic problem is outside the scope of this paper but will be addressed in future papers in collaboration with KAUST Extreme Computing Research Center (ECRC).
\vspace{-0.5em}
\begin{ack}
\vspace{-0.8em}
	I would like to express my deepest appreciation to Prof. G. Turkiyyah and Dr. S. Zampini without whom this paper would not have been possible.
\end{ack}
\vspace{-0.8em}
\bibliography{ifacconf}             
                                                   






\end{document}